\documentclass[11pt]{article}
\usepackage[english]{babel}
\usepackage{amsmath,amsthm}
\usepackage{amsfonts}
\usepackage{amssymb}
\usepackage{graphics}
\usepackage{graphicx}
\usepackage{color}
\usepackage{makecell}
\usepackage{pstricks}
\usepackage{ulem}

\theoremstyle{definition}

\theoremstyle{remark}

\everymath{\displaystyle}
\numberwithin{equation}{section}
\numberwithin{table}{section}
\definecolor{light-gray}{gray}{0.8}
\date{}
\begin{document}
\title{\textbf{Optimal solution of the liquidation problem under execution risk}}%
\author{Lorella Fatone\\ \small{Dipartimento di Matematica, Universit\`a di Camerino} \\ \small{Via Madonna delle Carceri 9, 62032
Camerino (MC), Italy}\\ \small{Ph. n. +39 0737 402558, Fax n. +39 0737 632525, E-mail: lorella.fatone@unicam.it}\\
Francesca Mariani \\ \small{Dipartimento di Scienze Economiche e Sociali, Universit\`a Politecnica delle Marche}\\ \small{Piazza Martelli 8, 60121 Ancona, Italy} \\\small{Ph. n. +39 071 2207243, Fax n. +39 071 2207102, E-mail: f.mariani@univpm.it} }
\maketitle
\begin{abstract}
We consider an investor that {trades continuously} and wants to liquidate {an initial asset position} within a prescribed time interval. During the execution of the liquidation order the investor  is subject to execution risk. 
We study the problem of finding the optimal liquidation strategy adopted by the investor in order to  maximize the expected revenue resulting from the liquidation.  We present a mathematical model of the liquidation problem that {extends} the model of  Almgren and Chriss (Almgren, R., Chriss, N.,  Optimal execution of portfolio transactions, Journal of Risk, 2000) to { include execution risk.}  The liquidation problem is modeled as a linear quadratic stochastic optimal control problem with finite horizon and, under some hypotheses, is solved explicitly. 
\vskip+0.3cm

\noindent{\bf Keywords:}
liquidation problem, stochastic optimal control, execution risk, Hamilton Jacobi Bellman equation\\
{\bf AMS Subject Classifications:} 93E20, 60H10, 49L20\\
{\bf JEL Codes:} C0, C61
\end{abstract}
\section{Introduction}\label{sec1}
The liquidation problem is the problem of finding the optimal strategy adopted by {an investor} in order to liquidate his position on a risky asset within a prescribed time interval, called liquidation interval. The liquidation problem is widely studied in mathematical finance (see, among others,  Almgren and Chriss, 2000, Almgren, 2003, Almgren, 2012, Ankirchner et al., 2016, Fatone et al., 2014,  Frei and Westray, 2015, Gu\'eant and Lehalle, 2015, Lorenz and Schied, 2012, Schied et al., 2010,  Tse et al., 2013). The mathematical model of the liquidation problem studied in these papers assumes that the execution of the liquidation order influences the asset share price inducing a slippage between the expected asset price and the actual price at which the trade is executed. In the financial markets usually this is the case when the liquidation order is a market order of large size. The earliest model of the liquidation problem has been introduced by Almgren and Chriss (2000), this is a discrete time model. Continuous time versions of the Almgren and Chriss model are developed by Almgren (2003),  Gatheral and Schied (2011) and Forsyth et al. (2012).  In these models the asset share price is the sum of an arithmetic Brownian motion and of a term that describes the impact of the {investor} trading activity. The utility function is the difference between the expected revenue resulting from the execution of the liquidation order and its variance. The liquidation problem is modeled as a mean variance optimization problem that is reduced to an elementary calculus of variations problem and solved. Several generalizations of the continuous time model introduced in Almgren (2003) have been developed. For example, Almgren (2012)  studies how liquidity affects the asset share price dynamics. In Fatone et al. (2014) the presence in the market of retail investors and its consequences on the execution of the liquidation order are considered. The retail investors are modeled as an homogeneous population of small investors whose behaviour is described by a mean field game. Gu\'eant and Lehalle (2015) assume the utility function to be a C.A.R.A. (Constant Absolute Risk Aversion) function and study the effects of limit order books on the execution of the order. In all these models the trading strategies are deterministic functions. Trading strategies that are stochastic processes have been considered in Schied et al. (2010), Ankirchner et al. (2016), Cheng et al. (2017), Bulthuis et al. (2017). In Schied et al. (2010) the trajectories of the trading  strategy are bounded and absolutely continuous functions of time defined in the liquidation interval. In Ankirchner et al. (2016) the effects of  trends in the asset share price on the execution of the liquidation order are studied and the trading rate is  modeled as a square integrable stochastic process of time.  In both papers a liquidation condition is imposed to the admissible trading strategies to require that at the end of the liquidation interval the initial asset share position is sold with probability one. The asset share price dynamic equation of  Schied et al. (2010)  and Ankirchner et al. (2016)  is the same used by Almgren (2003) and the liquidation problem is modeled as a stochastic optimal control problem. Under some hypotheses the value functions of the control problems are determined as solutions of the corresponding Hamilton Jacobi Bellman equations and of their auxiliary conditions (i.e. an initial condition in Schied et al., 2010, and a final condition in Ankirchner et al., 2016). The auxiliary condition used in Schied et al. (2010) and Ankirchner et al. (2016) is known in aeronautical engineering as fuel condition (see Bather and Chernoff, 1967, and, in the financial context, Schied et al. (2013)). The fuel condition of aeronautical engineering is a final condition that guarantees that no fuel is left unused at the end of the mission planned. In the liquidation problem the same condition guarantees that at the end of the liquidation interval the investor has completed (with probability one) the sale of the asset shares initially held. Trading strategies that are diffusion processes are considered in Cheng et al. (2017) to model the order fill uncertainty. The liquidation problem is solved in two different settings: in the first one the magnitude of order fill uncertainty is a prescribed positive constant parameter  independent from the trading strategy, in the second one the magnitude of order fill uncertainty is a linear function of the optimal trading rate. In the first setting the optimal trading strategy is found explicitly in terms of elementary functions without any constraints.  Otherwise, when the magnitude of uncertainty is a linear function of the trading rate, the system of Riccati equations associated to the Hamilton Jacobi Belmann equation is solvable under some strong assumptions on the parameters of the problem and the solution, when there exists, cannot be expressed in terms of elementary functions. More recently Bulthuis et al. (2017) have extended the model of Cheng et al. (2017) to include the uncertainty of limit order fills. The model is enriched by the addition of constraints to bound the trading rate of limit and market orders and of  a ``trade director'' to penalize trading strategies made simultaneously by  buy side market and sell limit orders. A further extension of the model of Cheng et al. (2017) is done by Cheng et al. (2019) in the case of constant uncertainty. The new model  adds to the old model a dynamic risk adjustment of the liquidation strategy. The risk adjustment is taken into account adding to the profit and loss function a quadratic term penalizing the strategies whose rate is  far from a prescribed target value. \\In this paper we present a mathematical model of the liquidation problem for {an investor whose trading strategy is subject to execution risk.}  {From an economic standpoint the execution risk can be generated by private taste shocks (Sannikov and Skrzypacz, 2016) or investor beliefs (Kyle et al., 2017) or private information regarding the asset value and/or inventories (Du and Zhu, 2017) as well as by uncertainty in order fills (Cheng et al., 2017, Bulthuis et al., 2017). When placing a market order, an investor is guaranteed to execute the order as the next available price. The actual price at which the order is executed can differ from the price at which the order is placed. This difference in financial trading is called slippage. Therefore an investor that schedules a market order gives a priority to the certainty of execution over the certainty of the execution price. However, there are not guarantees that the placed order, especially if large, is executed immediately. The causes of the lag between the placement and the settlement of an order  can be many,  from the unavailability of requested asset volume to the size of the order.  Similarly, private taste shocks or beliefs can deviate the realized trading strategy of the investor from the originally scheduled trading strategy. Because of the phenomena generating this deviation are hardly predictable, we refer to all them indifferently as execution risk.
As well as on the trading strategy,  the  execution risk impacts on the asset share price dynamics generating an additional source of risk different to  market risk. In line with Sannikov and  Skrzypacz (2016), Cheng et al. (2017), Bulthuis et al. (2017)  we model the effect of execution risk 
on the trading strategy assuming that the trading strategy is an \^Ito diffusion process whose noise term characterizes the magnitude of execution risk. The drift of the trading strategy is the scheduled trading rate and influences the asset share price dynamics. Moreover the impact of execution risk on the asset  share price is taken into account adding to the asset share price dynamic equation a noise term driven by a Wiener process correlated to the trading strategy. Both the noise terms of trading strategy and asset  share price dynamic equations are assumed to be square root functions of the scheduled trading rate and of the time left to reach the end of liquidation interval. {Note that, because of execution risk, at the end of the liquidation interval the {investor} can have a residual {asset position to sell in order to complete the liquidation order}}. In this case the residual {asset position} must be sold at the final time. In order to penalize the trading that at the end of the liquidation interval has not completed the liquidation, we consider as utility function of the control problem the sum of the expected revenue resulting from the liquidation and of a term penalizing the trading strategies that at the end of the liquidation interval have residual amount of asset shares left unsold. The asset share price dynamic equation of the model presented here is that of Almgren  and Chriss (2000) except for the temporary impact term that is proportional to the scheduled trading rate instead of the (actual) trading rate. The liquidation problem consists in finding the drift of the trading strategy (the scheduled trading rate) that maximizes the utility function.  The liquidation problem is formulated as a linear quadratic stochastic optimal control problem that has the {trading strategy} as state variable and the scheduled trading rate as control variable. To solve the model we use the completion of squares method to derive the Hamilton Jacobi Bellman equation and the optimal feedback control.  Explicit formulae of the optimal scheduled trading rate and of the corresponding optimal {trading strategy}  are found. The optimal scheduled trading rate of the model considered is determined and its dependence from the model parameters is studied. Unlike Cheng et al. (2017), the solution found has the advantage to be expressed in terms of elementary functions  and to be defined always independently from the values of the parameters appearing in the model. \\This paper is structured as follows. In Section 2 we formulate the liquidation problem. In Section 3, under some hypotheses on the form of execution risk, {we  solve the model introduced in Section 2.} In Section 4 we discuss some case studies that illustrate the behaviour of  model presented  in Section 2. Finally in Section 5 some conclusions are drawn.

\section{The model}\label{sec2}
We consider an investor that wants to liquidate within a fixed time interval, called liquidation interval, a prescribed number of shares of a risky asset traded in the financial market. Let $\mathbb R$ be the set of real numbers, $\mathbb R_+$ be the set of real positive numbers and $T,$ $Y\in\mathbb R$ be positive numbers. We denote by $[0,T]$ the liquidation interval and by $Y$  the initial amount of asset shares that {must be sold} within the time interval $[0,T].$ Let $y(t)$ be the trading strategy, i.e. the number of asset shares held by the investor at time $t\in[0,T],$ and $v(t,y(t)):[0,T]\times\mathbb{R}\rightarrow\mathbb{R}$ be the scheduled trading rate associated to the trading strategy $y(t),$ $t\in[0,T].$ To keep the notation simple, in the rest of paper the dependence of $v$ from $y$ is omitted and we use  the shorthand notation $v(t)$ to denote $v(t,y(t)),$ $t\in[0,T].$ The scheduled trading rate $v(t)$ is the rate with which the investor schedules to sell the asset shares at time $t,$ $t\in[0,T].$ Because of execution risk, the (realized) trading strategy $y(t),$ $t\in[0,T],$ satisfies the stochastic differential equation:
\begin{eqnarray}\label{1.1}
dy(t)=-v(t) \ dt+\phi(t,v(t))\ dW(t), \quad t\in[0,T],
\end{eqnarray}
where $\phi$ is a real function, $W(t),$ $t\in[0,T],$ is a standard Wiener process. The function $\phi$ characterizes the magnitude of execution risk. 
{The trading strategy of the investor in $[0,T]$ changes as consequence of the desire of the investor to buy or sell (modeled in (\ref{1.1}) by the term $vdt$) and as effect of execution risk 
(modeled in (\ref{1.1}) by the term $\phi dW$).}  {The presence in the trading strategy dynamics (\ref{1.1}) of the diffusion term $\phi dW$ has many possible explanations (Sannikov, Y., Skrzypacz,  2016, Cheng et al., 2017). For example, if the investor is a broker executing the liquidation order on the behalf of their clients, the diffusion term $\phi dW$ can model the shocks generated by the random orders of their clients. More generally, the diffusion term $\phi dW$ can model investor belief shocks (Kyle et al., 2017) or uncertainty in the order fills (Cheng et al., 2017 and Bulthuis et al., 2017)}. Equation (\ref{1.1}) is equipped with the initial condition:
\begin{eqnarray}
y(0)=Y.\label{1.2}
\end{eqnarray}
Equation (\ref{1.1}) is the state equation of the model of the liquidation problem studied in this paper, the initial condition (\ref{1.2}) assigns the {amount of asset shares} that must be sold within the liquidation interval $[0,T].$ The function $v$ is the unknown control variable of the liquidation problem. \\
Let $t\in[0,T],$ we denote by $S^0(t)$ the \textit{market} price of the asset share at time $t,$ and by $S(t)$  the  corresponding \textit{execution} price  (see Forsyth et al., 2012) at time $t,$  that is the price realized after the sale. We assume that $S(t),$ $t\in(0,T],$ is a stochastic process defined by the following equations:
\begin{align}\label{1.4}
&S(t) =S^0(t)+\kappa (H(t)-v(t)), &t\in(0,T],\\
&S^0(t)=S^0_0+\mu t + \gamma(y(t)-Y)+\psi B(t),& t\in(0,T], \label{1.6}
\end{align}
where $H(t)=\int_0^t\chi(s,v(s))dZ(s),$ $t\in[0,T],$ $\chi$ is a real function such that $\chi(t,0)=0,$ $t\in[0,T],$ $\mu\in\mathbb R,$ $\kappa,\gamma>0$ and $S_0^0>0$ are real constants, $B$ and $Z$ are Wiener processes. Note that the prices  $S^0(t),$ and $S(t),$ $t\in[0,T],$ solutions, respectively, of (\ref{1.4}),  (\ref{1.6}), are negative with positive probability. Usually this is an undesirable property since most of the time negative asset share prices are not allowed. However when $S_0^0>0$ and $\mu>0$ are large enough and sufficiently small values of $T$ are considered the event ``negative asset share prices'' has small probability and can be tolerated, as done in Almgren (2000, 2012),  Ankirchner et al. (2016), Fatone et al. (2014), Gu\'eant and Lehalle (2015), Cheng et al. (2017). \\ The terms $\gamma(y(t)-Y)$ and $\kappa(H(t)-v(t))$ are, respectively, the \textit{permanent} and  \textit{temporary impact} factors. \\The stochastic process $S^0(t),$ $t\in[0,T],$ describes the market price (see Cheng at al., 2017 and Forsyth et al., 2012), i.e. the asset share price that is not influenced from the temporary impact, and is defined by equation (\ref{1.6}), where $\psi>0$ is a real constant. The difference between the execution and the market prices is the so called slippage (see Almgren and Chriss, 2000).\\ We assume that:
\begin{eqnarray}
\mathbb E(dB(t), dW(t))=\mathbb E(dB(t), dZ(t))=0,\text{ and }\mathbb E(dZ(t),dW(t))=\rho dt,
\end{eqnarray} 
i.e. the Wiener processes $B(t),$ $W(t)$ and $B(t),$ $Z(t),$ $t\in[0,T],$ are uncorrelated and the Wiener processes $Z(t),$ $W(t),$ $t\in[0,T],$ are correlated with constant correlation coefficient $\rho\in[-1,1].$  \\ Equations (\ref{1.4}), (\ref{1.6}) extend the asset share price dynamic equation used by Almgren and Chriss and Chriss  (2000) to the case where the trading strategy is subject to execution risk. In (\ref{1.4}) the drift coefficient $v(t),$ $t\in[0,T],$ of the state equation (\ref{1.1}) replaces the time derivative of the {trading strategy} used in Almgren (2003). These two terms coincide when in (\ref{1.1}) we choose  $\phi\equiv 0.$ \\Let us justify the choice made in (\ref{1.4}) of using the term $\kappa v$ and the term $\kappa H$ to model the temporary impact factor, respectively,  {of the trading strategy} and of execution risk on the asset share price. First of all it must be said that when $\phi\neq 0$ the trajectories of the diffusion process (\ref{1.1}) are not differentiable, therefore it is not possible to consider their time derivative as done in Almgren and Chriss (2000). Second, it must be noted that,  in absence of the noise term $\phi dW(t)$ in (\ref{1.1}), the scheduled trading rate $v(t),$ $t\in[0,T],$ determines the trading strategy and, as consequence, affects the asset share price dynamics. Otherwise, when $\phi\neq0,$ in real markets, where the prices are the result of auctions, the desired rate of sale $v$ influences the asset share price dynamics even when, due to the unexpected  { circumstances} (modeled in (\ref{1.1}) with the term $\phi dW$), the desired rate of sale does not determine completely the {trading strategy dynamics}.  Choosing the scheduled trading rate $v,$ the {investor} chooses his desirable amount of asset shares to sell, this choice influences the asset share price dynamics and, consequently, the trading strategies of the {other} investors trading in the same asset. Differently from the term $\kappa v,$ the term $\kappa H$ affects directly the noise of the execution price $S$ and represents the additional noise generated by the influence of the scheduled order on the market price. This additional noise can be due, for example, to the behaviour, assumed by the other investors trading in the same asset, as consequence of the placement of the liquidation order.\\
The expected revenue resulting from the liquidation at time $T$ is given by:
\begin{eqnarray}\label{1.8}
\mathbb{E}\left[-\int_0^{T} (S(t)-S_0^0) dy(t)+(S(T)-S_0^0)y(T)\right],
\end{eqnarray}  
where $\mathbb E[\cdot]$ denotes the expected value of $\cdot.$ In (\ref{1.8}) the term 
\begin{eqnarray}\label{1.8b}
R=\mathbb{E}\left[(S(T)-S_0^0)y(T)\right]=\mathbb E\left[(S^0(T)-S^0_0)y(T)\right]+\mathbb E\left[-\kappa v(T)y(T)\right]
\end{eqnarray}
represents the expected revenue resulting from the liquidation at the final time $t=T.$ Since the trading strategy of the investor is subject to random noise (see Equation (\ref{1.1})) it is possible that at the end of the liquidation interval there is a residual amount of asset shares $y(T)$ to sell or buy, this adds to the expected revenue resulting from the liquidation at the market price $S^0(T)$ ($\mathbb E\left[(S^0(T)-S^0_0)y(T)\right]$) an extra cost, due to the risk of trading at the execution price $S(T)$ that is more unfavorable than the market price $S^0(T)$ (see Cheng at al., 2017). This extra cost is given by the term $\mathbb E\left[-\kappa v(T)y(T)\right].$ \\
In line with Cheng at al. (2017) we assume that:
 \begin{eqnarray}\label{1.8c}
 \mathbb E\left[-\kappa v(T)y(T)\right]=\mathbb E\left[-\lambda y^2(T)\right],
 \end{eqnarray}
 where $\lambda>0$ is a real constant.\\
For $t\in[0,T]$ given $\mathcal M_{[t,T]}$ be the set of the real-valued absolutely continuous and adapted processes in $[t,T]$, we define the set of admissible controls as the set of square integrable processes, that is: 
\begin{eqnarray}\label{opt_set}
\mathcal A_t=\left\{g\in \mathcal M_{[t,T]} \ : \ \int_t^T\mathbb E[g^2(t)]dt<+\infty \right\}.
\end{eqnarray}
The  liquidation problem is formulated as the following  linear quadratic stochastic optimal control problem:
\begin{eqnarray}\label{1.9}
\displaystyle\max_{v\in\mathcal{A}_0}\mathbb{E}\left[-\int_0^{T} (S(t)-S_0^0) dy(t)+(S^0(T)-S_0^0)y(T)-\lambda y^2(T)\right],
\end{eqnarray}
subject to the constraints (\ref{1.1}), (\ref{1.2}). \\
The penalization term $\mathbb E\left[-\lambda y^2(T)\right]$ in in (\ref{1.9})  measures the cost for selling  at time $T$ the residual amount of asset shares $y(T)$ at the execution price $S(T)$ instead of the market price $S^0(T).$  In line to what done by Karatzas et al. (2000) for the finite-fuel control problem and by Cheng et al. (2017) and Bulthuis et al. (2017)  for the liquidation problem, we consider a quadratic penalization term. It is worthing to note that as $\lambda\rightarrow +\infty$ the cost of selling at the end of the liquidation interval goes to infinity, i.e. the liquidation at time $T$ is not allowed and the final condition $y(T)=0$ is enforced. The condition $y(T)=0$ is the well known \textit{finite fuel constraint} introduced by Bene\v{s} et al. (1980) and further developed by Karatzas (1985).\\
When $\phi_0=\chi_0=0$ (i.e there is no execution risk) and $\lambda\rightarrow +\infty$ (i.e. the liquidation is completed at $T$ with probability one) problem (\ref{1.9}), (\ref{1.1}), (\ref{1.2})  reduces to the optimal execution problem solved by Almgren (2003), therefore the optimal trading strategy, solution of problem (\ref{1.9}), (\ref{1.1}), (\ref{1.2})  when $\mu=0,$ is the Volume Weighted Average Price (VWAP) strategy consisting in selling in each time interval an amount proportional to the predicted volume for that interval (Almgren, 2003).

\section{The solution}
In this section we solve problem  (\ref{1.9}), (\ref{1.1}), (\ref{1.2})  assuming:
\begin{align}
\phi(t,v)=&\phi_0\sqrt{(T-t)v}, \ t\in[0,T],\ v\in\mathbb R_+,\label{phi}\\
\chi(t,v)=&\chi_0\sqrt{(T-t)v},\ t\in[0,T],\ v\in\mathbb R_+,\label{chi}
\end{align} 
where $\phi_0,\chi_0>0.$ \\
As  already said in the Introduction and in Section 1, the diffusion terms $\phi(t,v),$ and $\chi(t,v),$ $t\in[0,T],$ $v\in\mathbb R_+,$ measure the magnitude of execution risk, respectively, in the investor trading strategy and in the asset share price dynamics. Whether due to uncertainty in  the order fills or to private taste shocks or beliefs, execution risk generates a deviation of the realized from the scheduled trading strategy and introduces an additional source of risk in the asset share price dynamics. In general, the larger the urgency to complete the liquidation order, and, as a consequence, the size of the residual asset position, the larger is the magnitude of execution risk  
(see Sannikov and Skrzypacz, 2016). With the choices (\ref{phi}), (\ref{chi}) we assume that  at time $t$ the magnitude of  execution risk affecting the trading strategy and the asset share price dynamics is proportional to the residual asset share position at time $t,$ that is  roughly of order $v(t)(T-t).$  Notice that  choices (\ref{phi}), (\ref{chi})   allow to obtain explicit solution of problem  (\ref{1.9}), (\ref{1.1}), (\ref{1.2})  expressed in terms of elementary functions without imposing any constraints on the parameters of the model.
\medskip

\noindent \textbf{Proposition 3.1}\\
Given $v\in\mathcal A_0,$ $S$ solution of (\ref{1.4}), (\ref{1.6}) and $y$ solution of (\ref{1.1}), (\ref{1.2}), the expected revenue in (\ref{1.9}) can be rewritten as follows:
\begin{eqnarray}\label{rev}
R=\mathbb{E}\left[-\lambda y^2(T)+\frac{\gamma}{2}(y^2(T)-Y^2)+\int_0^{T}\left(\mu y(t)+\left(\frac{\gamma}{2}\phi_0^2+\kappa\rho\chi_0\phi_0\right)(T-t)v(t)-\kappa v^2(t)\right)dt\right].
\end{eqnarray}
\textit{Proof.}\\
By   (\ref{1.1}), (\ref{1.2})  and (\ref{1.4}), (\ref{1.6}) we have:
\begin{align}\label{zero}
-\int_0^{T} (S(t)-S_0^0) dy(t)=&-\int_0^T(\mu t+\gamma(y(t)-Y)+\psi B(t)-\kappa v(t)+\kappa H(t))dy(t)\nonumber\\
=&-\mu\int_0^T t dy(t)-\gamma\int_0^T y(t)dy(t)+\gamma Y\int_0^T dy(t)-\psi\int_0^T B(t)dy(t)\nonumber\\
&+\kappa\int_0^T v(t) dy(t)-\kappa\int_0^TH(t)dy(t).
\end{align}
 Since:
\begin{align}
&y(t)dy(t)=\frac{1}{2}d(y^2(t))-\frac{1}{2}\phi_0^2(T-t)v(t)dt,& t\in[0,T],\label{uno}\\
&B(t)y(t)=d(B(t)y(t))-y(t)dB(t),& t\in[0,T],\label{due}\\
&H(t)dy(t)=d(H(t)y(t))-y(t)dH(t)-\rho\chi_0\phi_0(T-t)v(t)dt,& t\in[0,T].\label{tre}
\end{align}
Substituting (\ref{uno}), (\ref{due}), (\ref{tre}) into (\ref{zero})  we have:
\begin{align}\label{treb}
-\int_0^{T} (S(t)-S_0^0) dy(t)+&(S^0(T)-S_0^0)y(T)=-\kappa H(T)y(T)+\frac{\gamma}{2}(y^2(T)-Y^2)\nonumber\\
&+\int_0^T\left(\mu y(t)+\left(\frac{\gamma}{2}\phi_0^2+\kappa\rho\chi_0\phi_0\right)(T-t)v(t)-\kappa v^2(t)\right)dt\nonumber\\
&+\psi\int_0^T y(t)dB(t)+\kappa\phi_0\int_0^T\sqrt{(T-t)v^{3}(t)}dW(t)+\kappa\int_0^Ty(t)dH(t).
\end{align}
By the assumption $v\in\mathcal A_0,$  by the Jensen inequality and by (\ref{1.1}) there exists a real constant $K>0$ such that
$$ \sup_{t\in[0,T]}y^2(t)\leq  K\left(1+\int_0^T v^2(s)ds+\sup_{t\in[0,T]}\left(\int_0^t \phi_0^2(T-s)v(s)ds\right)^2\right)<+\infty,\ t\in[0,T],$$
applying the Burkholder-Davis-Gundy inequality there exists constants $K',K''>0$ such that
$$\mathbb E\left[\int_0^Ty^2(t)dt\right]\leq K'\mathbb E\left[\sup_{t\in[0,T]} y^2(t)\right]\leq K''\int_0^T\mathbb E\left[
1+\int_0^T(v(s)^2+\phi_0^2(T-s)v(s))ds\right]<+\infty.$$
Then we have:
 \begin{eqnarray}\label{quattro}
 \mathbb E\left[\int_0^T y(t)dB(t)\right]=0.
 \end{eqnarray}
By  the assumption $v\in\mathcal A_0$ we have
$\mathbb E\left[\int_0^T(T-t)v^{3}(t)dt\right]\leq T\mathbb E\left[\int_0^T v^3(t)dt\right]<\infty$
then 
\begin{eqnarray}\label{cinque}
\mathbb E\left[\int_0^T\sqrt{(T-t)v^{3}(t)}dW(t)\right]=0.
\end{eqnarray}
Moreover from
$\mathbb E\left[\int_0^T(T-t)^2v^2(t)dt\right]\leq T^2\mathbb E\left[\int_0^T v^2(t)dt\right]<+\infty$
it follows that the stochastic process $H(t),$ $t\in[0,T],$ is a martingale and  
\begin{eqnarray}\label{sei}
\mathbb E\left[\int_0^T y(t)dH(t)\right]=0.
\end{eqnarray}
Finally substituting (\ref{quattro}), (\ref{cinque}), (\ref{sei}) into (\ref{treb}) we obtain (\ref{rev}). This concludes the proof.
\bigskip

\noindent{\textbf{Proposition 3.2}}\\
The value function of stochastic optimal control problem  (\ref{1.9}), (\ref{1.1}), (\ref{1.2}) satisfies the following Hamilton Jacobi Bellmann equation:
\begin{eqnarray}\label{HJB}
\frac{\partial V(t,y)}{\partial t}+\frac{1}{4\kappa}\left(\frac{\phi_0^2}{2}(T-t)\frac{\partial^2 V(t,y)}{\partial y^2}+\left(\frac{\gamma}{2}\phi_0^2+\kappa\rho\chi_0\phi_0\right)(T-t)-\left(\frac{\partial V(t,y)}{\partial y}+\gamma y\right)\right)^2+\mu y=0
\end{eqnarray}
with final condition:
\begin{eqnarray}\label{sc}
V(T,y)=-\lambda y^2.
\end{eqnarray}
The optimal scheduled {trading rate} $v^*(t),$ $t\in[0,T],$ solution of problem (\ref{1.9}), (\ref{1.1}), (\ref{1.2}) has the state-feedback expression:
\begin{eqnarray}\label{rate_opt}
v^*(t,y)=&\frac{y(t)}{T-t+\alpha}-\frac{1}{4\kappa}(\mu+B)\left(T-t+\alpha-\frac{\alpha^2}{T-t+\alpha}\right)+\frac{1}{2\kappa}B(T-t)\nonumber\\&+\frac{\alpha}{2\kappa}\left(\frac{B(T-t)}{T-t+\alpha}-\frac{\kappa\phi_0^2}{T-t+\alpha}\ln\left(\frac{T-t+\alpha}{\alpha}\right)\right), \ t\in[0,T], \ y\in\mathbb R,
\end{eqnarray}  
where $\alpha=\frac{2\kappa}{2\lambda-\gamma}>0$ and $B=\frac{\gamma}{2}\phi_0^2+\kappa\rho\chi_0\phi_0.$\\
\textit{Proof} \\ We use the \textit{completion of squares} method (see Brokett, 1970). Using (\ref{rev}) the liquidation problem becomes:
\begin{eqnarray}\label{1.9b}
\displaystyle\max_{v\in\mathcal{A}_0}\mathbb E\left[-\lambda y^2(T)+\frac{\gamma}{2}(y^2(T)-Y^2)+\int_0^{T}\left(\mu y(t)+\left(\frac{\gamma}{2}\phi_0^2+\kappa\rho\chi_0\phi_0\right)(T-t)v(t)-\kappa v^2(t)\right)dt\right]
\end{eqnarray}
subject to constraints (\ref{1.1}), (\ref{1.2}). 
The value  function associated to problem (\ref{1.9b}), (\ref{1.1}), (\ref{1.2}) is given by:
\begin{eqnarray}\label{3.1}
V(t,y)=&\displaystyle\max_{v\in\mathcal{A}_t}\mathbb{E}_t\left[-\lambda y^2(T)+\frac{\gamma}{2}(y^2(T)-y^2(t))+\int_t^{T}\left(\mu y(s)+\left(\frac{\gamma}{2}\phi_0^2+\kappa\rho\chi_0\phi_0\right)(T-s)v(s)\right.\right.\nonumber\\
&\left.\left.-\kappa v^2(s)\right)ds\right],\ t\in[0,T],
\end{eqnarray}
where the maximum is taken over the class of the trading strategies  solutions of (\ref{1.1}), (\ref{1.2}) whose scheduled trading rate belongs to  $\mathcal A_t.$ In (\ref{3.1}) $\mathbb E_t[\cdot]$ denotes the conditional expectation $\mathbb E[\cdot | y(t)=y],$ $t\in[0,T].$ \\Applying \^Ito formula to $y^2(t),$ $t\in[0,T],$ and using (\ref{1.1}), (\ref{1.2}) we have :
\begin{eqnarray}\label{3.2}
y^2(T)=y^2(t)+\int_t^T\left(\phi_0^2(T-s)v(s)-2y(s)v(s)\right)ds+2\int_t^T\phi_0\sqrt{(T-s)v(s)}y(s)dW(s), \ t\in[0,T],\
\end{eqnarray}
then the value function $V$ in (\ref{3.1}) reduces to:
\begin{eqnarray}\label{3.4}
V(t,y)=-\kappa\displaystyle\min_{v\in\mathcal{A}_t}\mathbb{E}_t&\left[\int_t^T\left(v^2(s)+\frac{2}{\alpha}\gamma y(s)v(s)-\frac{1}{\kappa}\left(\phi_0^2\gamma+\kappa\rho\chi_0\phi_0\right)(T-s)v(s)\right.\right.\nonumber\\&\left.\left.\quad-\frac{1}{\alpha}\phi_0^2(T-s)v(s)-\frac{\mu}{\kappa}y(s)\right)ds+\frac{\lambda}{\kappa}y^2\right],\ t\in[0,T],
\end{eqnarray}
where $\alpha=\frac{2\kappa}{2\lambda-\gamma}.$\\
Let:
\begin{align*}
f_1(t)=&-\frac{1}{2\kappa}(\mu+B)\left(T-t+\alpha-\frac{\alpha^2}{T-t+\alpha}\right)+\frac{\phi_0^2}{T-t+\alpha}(T-t)\nonumber\\
&+\alpha\left(\frac{1}{\kappa}B(T-t)+\frac{\phi_0^2}{T-t+\alpha}\right)\ln\left(\frac{T-t+\alpha}{\alpha}\right),\ t\in[0,T],\\
f_2(t)=&\frac{1}{T-t+\alpha}-\frac{1}{\alpha},\ t\in[0,T].
\end{align*}
we observe that $f_1(T)=f_2(T)=0,$ $f_2'(t)=1/(T-t+\alpha)^2$ and
\begin{eqnarray*}
f_1'(t)=\frac{1}{T-t+\alpha}\left(\frac{\phi_0^2}{T-t+\alpha}(T-t)-\frac{1}{\kappa}B(t-t)+ f_1(t)\right)-\mu,\ t\in[0,T].
\end{eqnarray*}
Applying  \^Ito formula to $f_1(t)y(t)$ and $f_2(t)y^2(t),$ $t\in[0,T],$  yields:
\begin{align}
0=&f_1(t)y(t)-\int_t^T \left(f_1(s)v(s)-f_1'(s)y(s)\right)ds+\int_t^T f_1(s)\phi_0 \sqrt{(T-s)v(s)}dW(s),\ t\in[0,T],\label{3.5}\\
0=&f_2(t)y^2(t)-\int_t^T\left(2f_2(s)v(s)y(s)-\frac{y^2(s)}{(T-s+\alpha)^2}-f_2(s)\phi_0^2(T-s) v(s)\right)ds\nonumber\\&+\int_t^T 2f_2(s)\phi_0\sqrt{(T-s)v(s)} dW(s), \ t\in[0,T].\label{3.6}
\end{align}
Since $f_1(t)$ and $f_2(t)$ are bounded in $[0,T]$ the stochastic integrals of $\displaystyle \int_t^T f_1(s)\phi_0 ((T-s)v(s))^{1/2}dW(s)$ and $\displaystyle \int_t^T 2f_2(s)\phi_0\sqrt{(T-s)v(s)} dW(s)$ has zero expectation (thought they are not necessarily  martingales) and from (\ref{3.4}), (\ref{3.5}), (\ref{3.6}) we have:
\begin{align}\label{3.7}
V(t,y)=-\kappa\displaystyle\min_{v\in\mathcal{A}_t}\mathbb{E}_t&\left[\int_t^T\left(v^2(s)+\frac{2}{\alpha}\gamma y(s)v(s)-\frac{1}{\kappa}\left(\phi_0^2\gamma+\kappa\rho\chi_0\phi_0\right)(T-s)v(s)-\frac{1}{\alpha}\phi_0^2(T-s)v(s)\right.\right.\nonumber\\&\left.\left.\quad-\frac{\mu}{\kappa}y(s)\right)ds+\frac{\lambda}{\kappa}y^2\right]+\mathbb E_t\left[f_1(T)y(T)+f_2(T)y^2(T)\right]\nonumber\\
=&-\kappa\displaystyle\min_{v\in\mathcal{A}_t}\mathbb{E}_t\left[\int_t^T\left(v^2(s)+\frac{2}{\alpha}\gamma y(s)v(s)-\frac{1}{\kappa}\left(\phi_0^2\gamma+\kappa\rho\chi_0\phi_0\right)(T-s)v(s))\right.\right.\nonumber\\&\left.\left.\quad-\frac{1}{\alpha}\phi_0^2(T-s)v(s)-\frac{\mu}{\kappa}y(s)+f_1'(s)y(s)-f_1(s)v(s)+\frac{y^2(s)}{(T-s+\alpha)^2}\right.\right.\nonumber\\&\left.\left.
-\frac{2}{T-s+\alpha}y(s)v(s)+\frac{2}{\alpha}y(s)v(s)+\frac{\phi_0^2(T-s)}{T-s+\alpha}v(s)-\frac{\phi_0^2(T-s)}{\alpha}v(s)\right)ds\right.\nonumber\\&\left.-f_1(t)y(t)+\left(\frac{1}{T-t+\alpha}-\frac{1}{\alpha}+\frac{\lambda}{\kappa}\right)y^2(t)\right],\ t\in[0,T].
\end{align} 
Now, adding and subtracting to (\ref{3.7}) the term:
$\frac{1}{4\kappa^2}\int_t^T\left(\left(-\frac{\kappa\phi_0^2}{T-s+\alpha}+B\right)(T-s)-f_1(s)\right)^2ds, \ t\in[0,T],$
we obtain:
\begin{align}\label{3.10}
V(t,y)= &-\kappa\min_{v\in\mathcal{A}_t}\mathbb E_t\left[\int_0^T\left(v(s)-\frac{1}{T-s}y(s)-\frac{1}{4\kappa}\left(\frac{\gamma}{2}\phi_0^2+\rho\kappa\chi_0\phi_0-\mu\right)(T-s)\right)^2ds\right]\nonumber\\ &+c(t)-\kappa f_1(t)y-\left(\frac{\kappa}{T-t+\alpha}+\frac{\gamma}{2}\right)y^2,
\ t\in[0,T],\ y\in\mathbb R,
\end{align}
where $c'(t)=-\frac{1}{4\kappa}\left(\left(-\frac{\kappa\phi_0^2}{T-t+\alpha}+B\right)(T-t)-f_1(t)\right)^2,$ $t\in[0,T].$\\
By straightforward computations it is easy to verify that  the maximum in (\ref{3.10}) is attained at
 $v=v^*$ where $v^*$ is given by (\ref{rate_opt}) and
 \begin{eqnarray}
V(t,y)= a(t)y^2+b(t)y+c(t), \ t\in[0,T],\ y\in\mathbb R,
\end{eqnarray}
where:
\begin{align}
a(t)=&-\frac{\gamma}{2}-\frac{\kappa}{T-t}, \ t\in[0,T], \label{3.11}\\
b(t)=&-\kappa f_1(t),\ t\in[0,T].\label{3.12}
\end{align}
Note that the functions $a(t),$ $b(t),$  $t\in[0,T],$ are solutions of the following system of Riccati equations:
\begin{align}
a'(t)=&-\frac{1}{\kappa}\left(a(t)+\frac{\gamma}{2}\right)^2,& t\in[0,T],\label{3.14}\\
b'(t)=&\frac{1}{\kappa}\left(a(t)+\frac{\gamma}{2}\right)\left(\phi_0^2(T-t)\left(a(t)+\frac{\gamma}{2}\right)+B(T-t)-b(t)\right)-\mu, &t\in[0,T]\label{3.15}
\end{align}
with final conditions: $a(T)=-\lambda,$ $b(T)=0.$ \\
Finally by straightforward computations it is easy to verify that  the value function $V$ satisfies the Hamilton Jacobi Bellmann equation (\ref{HJB}) with final condition (\ref{sc}). This concludes the proof.
\bigskip

\noindent\textbf{Corollary 3.1}\\
In the limit as $\lambda\rightarrow +\infty$ the optimal scheduled trading rate reduces to:
\begin{eqnarray}\label{opt}
v^*(t,y)=\frac{y(t)}{T-t}-\frac{1}{4\kappa}(\mu-B)(T-t), \ t\in[0,T], \ y\in\mathbb R.
\end{eqnarray}
\textit{Proof.} It follows easily taking the limit of (\ref{rate_opt}) as $\lambda\rightarrow +\infty$.
\bigskip

\noindent
Recall that when $\phi_0=0$ the optimal scheduled trading rate $v^*$ in (\ref{opt}) is the optimal trading rate found by Almgren (2003) under constant directional view about the asset price evolution (see Ankirchner et al. 2016). In the case where we have also zero drift ($\mu=0$) the optimal trading strategy $y^*$ is the VWAP strategy that consists to sell in each time interval an amount of asset shares proportional to the predicted volume for that interval (see Almgren, 2003). On the other hand, it is worth to note that, when $\phi_0\neq 0,$ the optimal scheduled trading rate in (\ref{opt}) is the optimal trading rate of Almgren (2003) for a modified asset price $S$ with drift given by $\tilde \mu=\mu-B=\mu-\frac{\gamma}{2}\phi_0^2-\kappa\rho\chi_0\phi_0.$  In other words, under execution risk the investor modifies his directional view about the future asset price growth rate passing from $\mu$ to $\mu-B.$ It should be noted that, when $\rho\geq 0$ (i.e. there is non negative correlation between trading strategy and asset share price dynamics)  the asset drift $\tilde\mu$ in presence of execution risk is smaller than the asset drift $\mu$ in absence of execution risk. Otherwise when $\rho<0$ (i.e. there is negative correlation between trading strategy and asset share price dynamics) we have $\tilde\mu>\mu.$ After all, it is legitimate to believe that asset price and trading strategy are positive correlated. In fact, when execution risk affects the trading strategy determining a decrease on the amount of asset shares sold with respect to the scheduled amount, we expect an increase in the asset share price. Therefore, assuming a positive correlation between asset share price and trading strategy,  we can conclude that the presence of execution risk changes the directional view of the investor regarding the future  price movement  causing him to expect a lower asset share return than in absence of execution risk.
\bigskip

\noindent \textbf{Proposition 3.3}\\
Let $y^*(t),$ $t\in[0,T],$ be the optimal trading strategy of problem (\ref{1.9b}), (\ref{1.1}), (\ref{1.2}) as $\lambda\rightarrow +\infty$ we have
 $\lim_{t\rightarrow T^-} y^*(t)=0$ a.s..\\
\textit{Proof.}\\
Here we follow Delyon and Hu (2006).
From Corollary 3.1 substituting $v^*,$ given by formula (\ref{opt}), into (\ref{1.1}) we obtain that  the optimal {trading strategy} $y^*(t),$ $t\in[0,T],$ associated to problem (\ref{1.9b}), (\ref{1.1}), (\ref{1.2}) with $\phi(t,v)=\phi_0((T-t)v)^{1/2}$ and  $\chi(t,v)=\chi_0((T-t)v)^{1/2},$ $t\in[0,T],$ $v\in\mathbb R_+,$  is solution of the following problem:
\begin{align}\label{3.17}
dy^*(t)=&-\left(\frac{y^*(t)}{T-t}-\frac{1}{4\kappa}(\mu-B)(T-t)\right)dt+\phi_0\sqrt{y^*(t)-\frac{1}{4\kappa}(\mu-B)(T-t)^2}dW(t), \ t\in[0,T],\\
y^*(0)=&Y.\label{3.18}
\end{align}
Let $\tilde y(t)=y^*(t)-\frac{1}{4\kappa}(\mu-B)(T-t)^2,$ $t\in[0,T],$ applying \^Ito's formula to $y^*(t),$ $t\in[0,T],$ it is easy to verify by straightforward computations that $\tilde y$ is solution of:
\begin{align}\label{3.19}
d\tilde y(t)=&-\left(\frac{\tilde y(t)}{T-t}-\frac{1}{2\kappa}(\mu-B)(T-t)\right)dt+\phi_0\sqrt{\tilde y(t)}dW(t), \ t\in[0,T],\\
\tilde y(0)=&\tilde y_0,\label{3.20}
\end{align}
where $\tilde y_0=Y-\frac{1}{4\kappa}(\mu-B)T^2.$\\
Applying \^Ito's formula to $\frac{\tilde y(t)}{T-t},$ $t\in[0,T],$ we deduce:
\begin{eqnarray}\label{3.21}
\frac{\tilde y(t)}{T-t}=\frac{\tilde y_0}{T}+\frac{1}{2\kappa}(\mu-B)t+\int_0^t\frac{\sqrt{\tilde y(s)}}{T-s} dW(s),\ t\in[0,T].
\end{eqnarray}
Since the stochastic process $\left\{\frac{\sqrt{\tilde y(t)}}{T-t}\right\}_{t\in[0,T]}$  is locally bounded a.s., then $M(t)=\int_0^t\frac{\sqrt{\tilde y(s)}}{T-s} dW(s),$ $t\in[0,T],$ is a  
martingale with quadratic variation:
\begin{eqnarray}\label{3.22}
\langle M \rangle(t)=\int_0^t \frac{\tilde y(s)}{(T-s)^2} ds,\ t\in[0,T].
\end{eqnarray}
Note that $\langle M \rangle(t)\rightarrow +\infty$ as $t\rightarrow T^-$ and there exists a constant $K>0$ such that $\langle M\rangle(t)\leq \frac{K}{T-t},$ $t\in[0,T].$\\
Applying Dambis–Dubins–Schwarz’s theorem (see Klebaner, 2012), we have that there exists a standard one-dimensional Brownian motion $\hat B$ such that:
\begin{eqnarray}\label{3.23}
M(t)=\hat B(\langle M\rangle(t)),\ t\in[0,T].
\end{eqnarray}
Substituting (\ref{3.23}) into (\ref{3.21}) we have:
\begin{eqnarray}\label{3.24}
\tilde y(t)=(T-t)\left(\frac{\tilde y_0}{T}+\frac{1}{2\kappa}(\mu-B)t+\phi_0\hat B(\langle M\rangle(t))\right),\ t\in[0,T].
\end{eqnarray}
Finally, since the limit of $t\hat B(1/t)$ as $t\rightarrow 0$ goes to zero by the Law of Large Numbers for Brownian motions, we have that:
\begin{eqnarray}\label{3.25}
 \lim_{t\rightarrow T^-} (T-t)\hat B(\langle M\rangle(t))=0 \ \text{a.s.},
\end{eqnarray}
and
\begin{eqnarray}\label{3.26}
\lim_{t\rightarrow T^-} y^*(t)=\lim_{t\rightarrow T^-} \left(\tilde y(t)+\frac{1}{4\kappa}(\mu-B)(T-t)^2\right)=0 \ \text{a.s.}.
\end{eqnarray}
This concludes the proof.
\bigskip

\noindent The process $\tilde y,$ solution of the stochastic differential equation (\ref{3.19}), is the diffusion process of  Deylon and Hu (2206) constructed by adding to the process $\hat y$, solution of $d\hat y(t)=\phi_0\sqrt{\hat y(t)},$ $t\in[0,T],$ the extra drift term $-\hat y(t)/(T-t)+1/2\kappa(\mu-B)(T-t).$ As $t\rightarrow T^-$  this last term becomes increasingly strong forcing the process $\tilde y$ to hit $0$ at $t = T$ a.s. (see Deylon and Hu, 2006, and Whitaker et al., 2016). When $B=\mu$ a popular discretization of the stochastic differential equation (\ref{3.19}) is the Modified Diffusion Bridge introduced by Durham and Gallant (2002). Notice that the process $\tilde y,$ solution of (\ref{3.19}), is absolutely continuous with respect to the conditioned process  $\hat y|0,$ that is the process $\hat y$ conditioned on hitting $0$ a.s. at $t = T.$ \\The processes $y^*(t),$ $\tilde y(t),$ $t\in[0,T],$ solutions of (\ref{3.17}), (\ref{3.18}) and of (\ref{3.19}), (\ref{3.20}), are Extend Cox Ingersoll Ross (ECIR) square root processes (Hull and White, 1990)  with reversion rate $-1/(T-t)$ and time dependent equilibrium levels given, respectively, by $\frac{1}{4\kappa}(\mu-B)(T-t)^2$ and $\frac{1}{2\kappa}(\mu-B)(T-t)^2.$ \\By straightforward computations we obtain that the expected value of $y^*(t),$ $t\in][0,T],$ is given by $\mathbb E(y^*(t))=\left(\frac{Y}{T}+\frac{1}{4\kappa}(\mu-B)t\right)(T-t),$ $t\in[0,T].$ When $\mu>B,$ i.e. when the asset growth rate $\mu$ is large enough or the execution risk parameter $\phi_0$ is small enough, the expected value of the optimal strategy is a concave function of time, this means that the investor on average liquidates the initial asset position more quickly over time.  This is the behaviour of an investor  believing that the asset price will rise in the future, and, as a consequence, postpones selling in time to take advantage of the asset price increase. Otherwise, when $\mu<B,$  i.e. when the asset growth rate $\mu$ is small enough or the execution risk parameter $\phi_0$ is large enough, the expected value of the optimal strategy is a convex function of time, this means that the investor on average liquidates the initial asset position more slowly over time.  This is the behaviour of an investor  believing that the asset price is likely to decrease in the future, and, as a consequence, sells more quickly at the beginning of the liquidation to avoid disadvantages of the asset price decrease. \\
Differently to Cheng et al. (2017), where the risk uncertainty affects the optimal trading strategy only in its diffusion term, in our model the risk uncertainty affects also the drift of the optimal trading strategy changing the directional view of the investor about the price movement. It is interesting to observe that when we choose $\gamma=2\kappa\rho\chi_0$ we have $B=0$ and, in this case, the drift of the optimal trading strategy $y^*,$ solution of (\ref{3.17}), (\ref{3.18}), does not depend on $\phi_0,$ this happens only if we choose $\rho>0,$ i.e. if we assume that asset share price and trading strategy are positive correlated.

\section{Case studies}
In this section we analyze the behaviour of the optimal trading strategy obtained in Proposition 3.2 in two case studies that differ for the order considered. Moreover we compare the optimal trading strategy obtained in Proposition 3.2  with the adaptive VWAP strategy (also called constant uncertainty trading strategy) of Cheng et al. (2017). The adaptive VWAP strategy is the solution of problem (\ref{1.9}), (\ref{1.1}), (\ref{1.2}) in the case where the execution risk parameters are given by $\chi_0=0$ and  $\phi=m_0,$ where $m_0$ is a real constant. When $\phi_0=\chi_0=0$ the optimal trading strategy obtained in Proposition 3.2 and the adaptive VWAP strategy of Cheng et al. (2017) coincide with the deterministic VWAP strategy of Almgren and Chriss (2000). For shortness, in the rest of Section we call the optimal trading strategy and rate obtained in Proposition 3.2 square root uncertainty trading strategy and rate. \\We simulate, with the explicit Euler method, the optimal trading strategy, solution of (\ref{1.9}), (\ref{1.1}), (\ref{1.2}), and the adaptive VWAP strategy of Cheng et al. (2017) using as parameters of simulation those used in Almgren and Chriss (2000) and Cheng et al. (2017). To guarantee a  fair comparison between the two models, across all simulations we generate the trajectories using the same Brownian motions. Specifically,  assuming the trading year made by $252$ trading days, we consider as time unit a trading day and we choose:  the initial asset share position to liquidate $Y=10^6,$ the liquidation interval of one day $T=1,$ the initial asset share price $S_0^0=50\$/share$, the permanent impact parameter $\gamma=2.5\times 10^{-7}\$/share^2$, the temporary impact parameter $\kappa=2.5\times 10^{-6}(\$/share^2)day$ and $\lambda=1000\kappa.$
\begin{figure}[hptb]
	\centerline{\includegraphics[height=7cm]{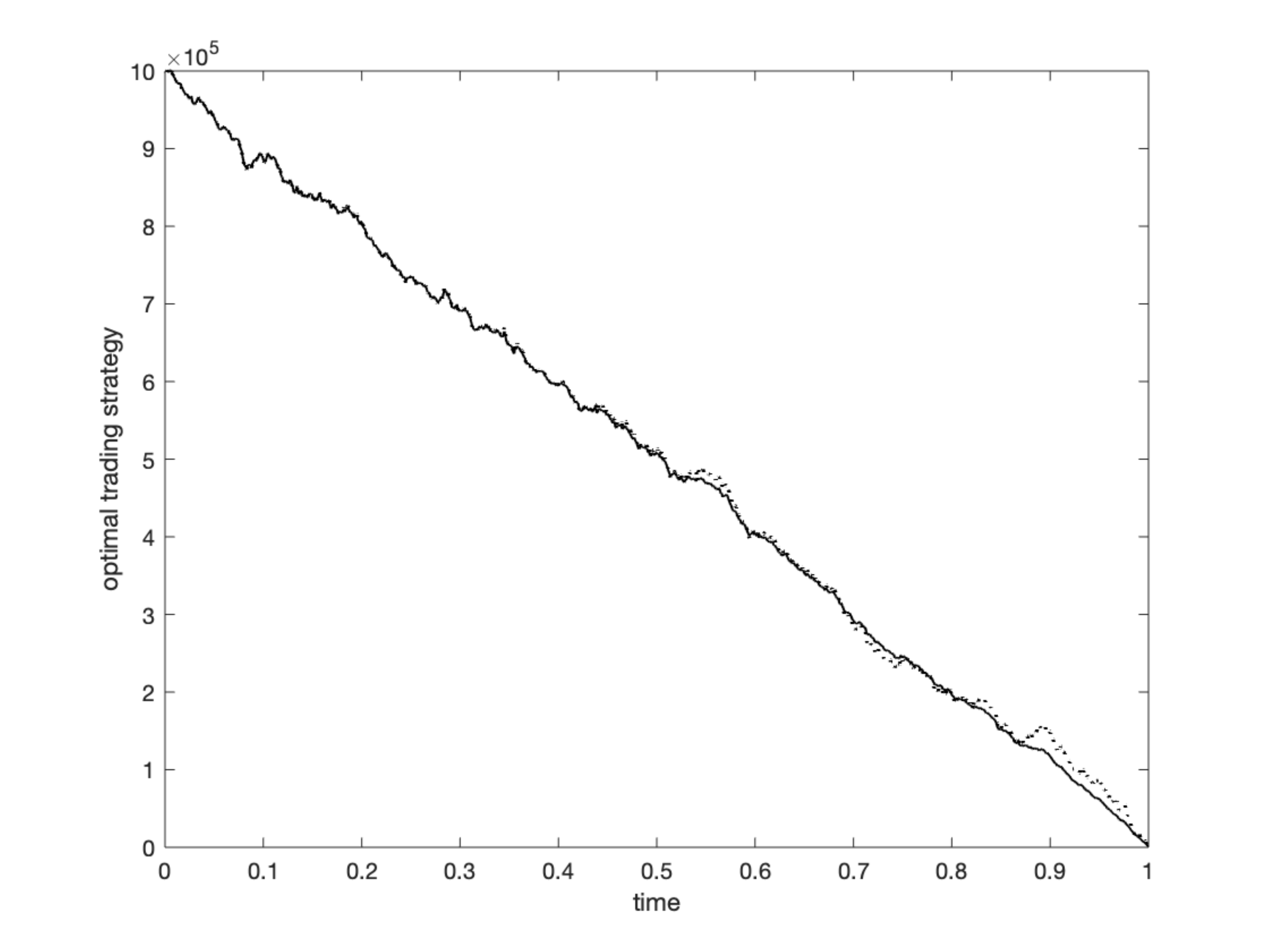}\includegraphics[height=7cm]{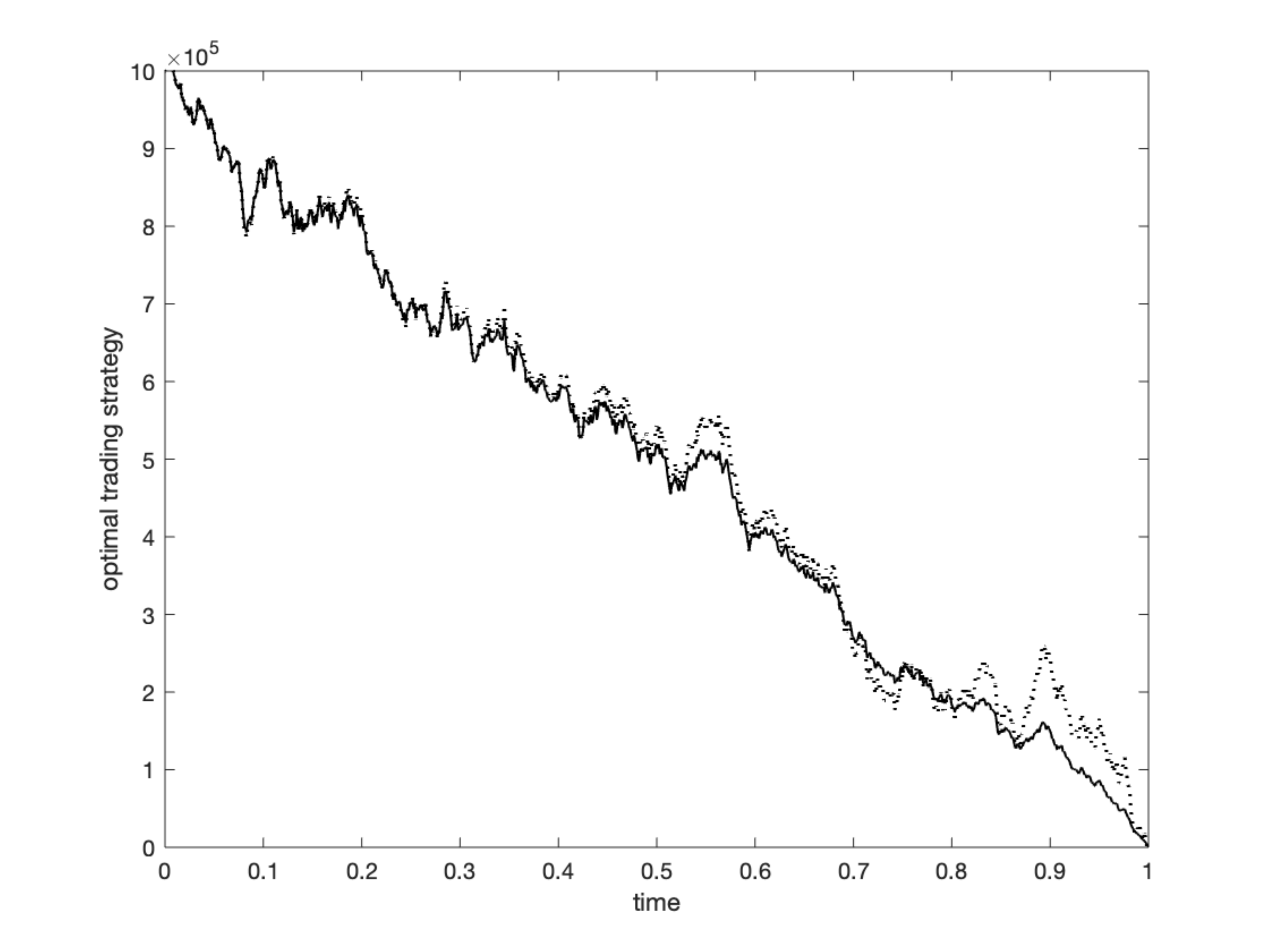}}
	\caption{Sample trajectories of optimal square root uncertainty optimal trading strategy (solid line) and of the optimal constant uncertainty optimal trading strategy (dotted line) obtained with $p_0=10\%$ (left panel) and $p_0=30\%$ (right panel).  }\label{fig1}
	\centerline{\includegraphics[height=7cm]{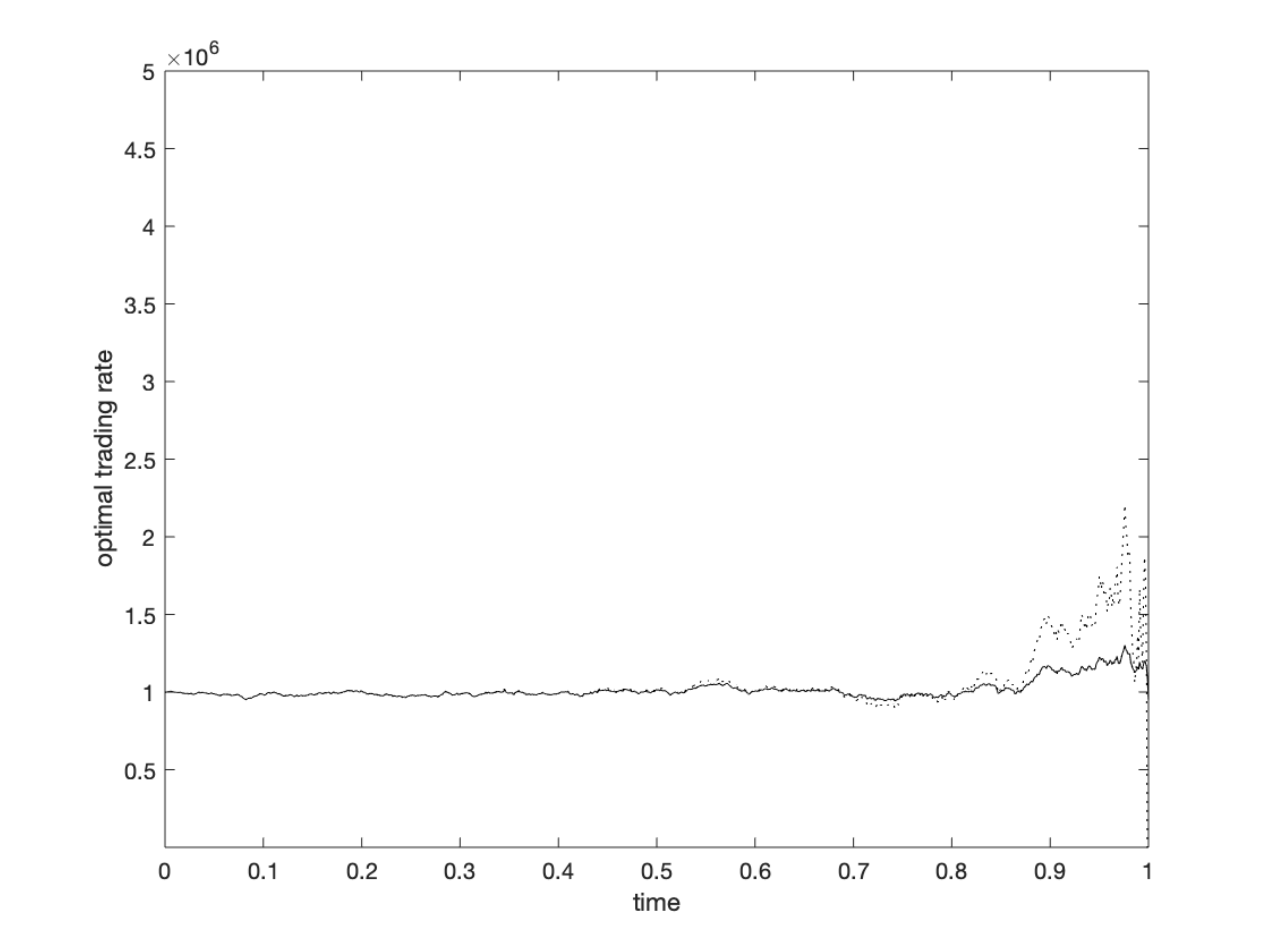}\includegraphics[height=7cm]{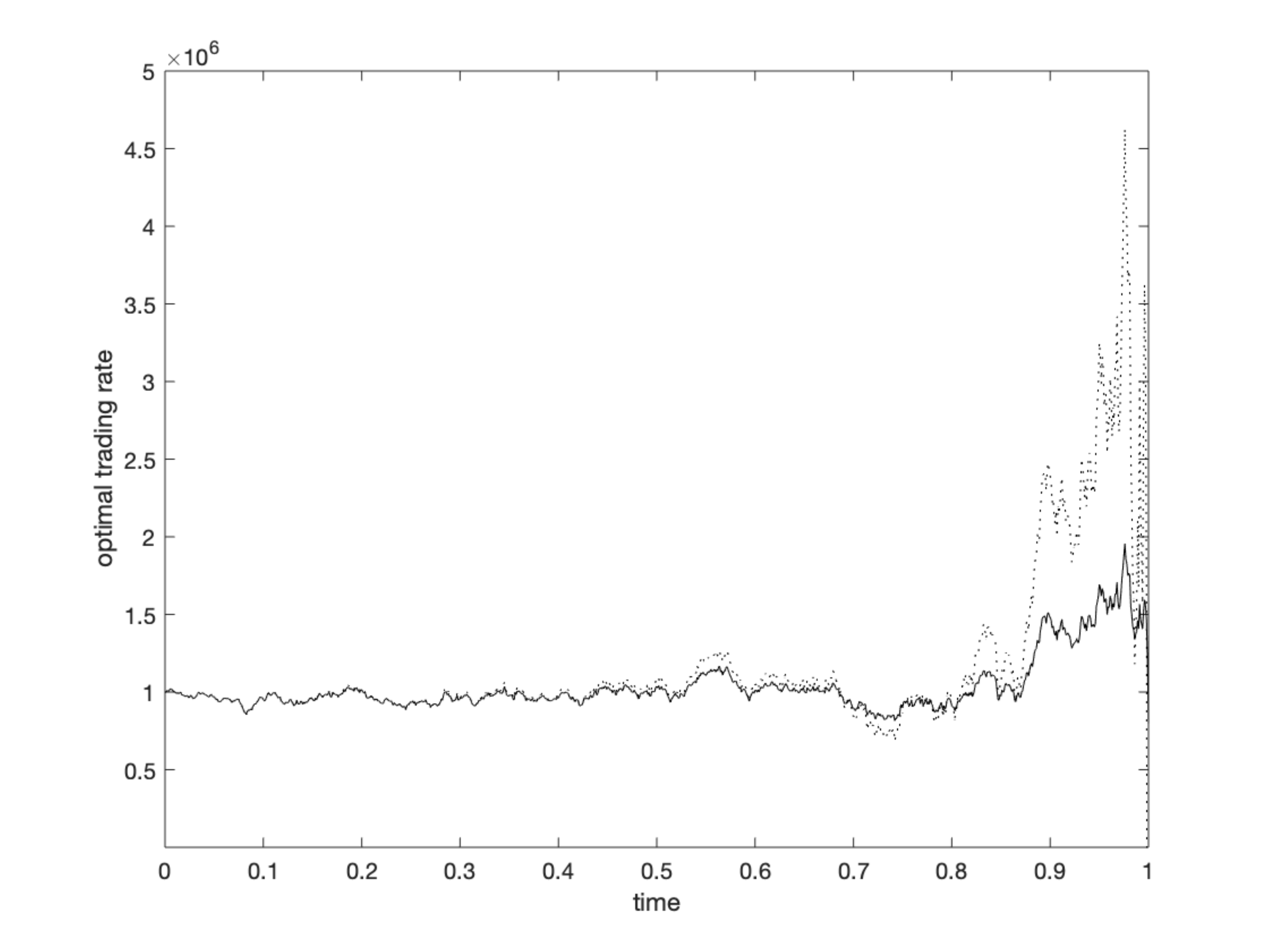}}
	\caption{Sample trajectories of optimal square root uncertainty optimal trading rate (solid line) and of the optimal constant uncertainty optimal trading rate (dotted line) obtained with $p_0=10\%$ (left panel) and $p_0=30\%$ (right panel).  }\label{fig2}
\end{figure}
Let  $p_0>0,$ in line with Cheng et al. (2017), we choose $m_0=p_0 Y,$ $\phi_0=p_0\sqrt{Y/T}share^{1/2}$ and $\chi_0=\sigma.$  With these choices the executed orders have on average $p_0$ deviation from the placed orders per day. In fact at each time $t\in[0,T]$ the constant $\phi_0$ multiplies $\sqrt{v(t)(T-t)},$ where $v(t)$  is roughly of order $Y/T$. Given $Y$ and $T$ the difficulty of liquidation increases when the execution risk parameter  $p_0$ increases. 
The aim of this section is to analyze the behaviour of the optimal trading strategy obtained in Proposition 3.2 when  ``easy'' and ``difficult'' orders are considered. We consider as ``difficult''  a liquidation order with small values of $p_0$ and we consider ``easy''  a liquidation order  with large values of $p_0.$ Specifically in the numerical experiments we choose $p_0=10\%$ for the ``easy''  order and $p_0=30\%$ for the ``difficult''  order. \\ In  Figure 1 we plot the sample trajectory of the optimal square root uncertainty trading strategy (solid line) and of the optimal constant uncertainty optimal trading strategy (dotted line) obtained with $p_0=10\%$ (left panel) and $p_0=30\%$ (right panel).  Looking at Figure 1 we observe that at the beginning of the liquidation interval the optimal strategies are very close each other and close to the VWAP strategy that corresponds to a linear reduction of holdings over the liquidation interval.  As time approaches to the liquidation horizon the optimal square root uncertainty trading strategy moves away from the optimal constant uncertainty strategy and, except for the time interval $[0.7,0.8],$ is under the optimal constant uncertainty strategy. This behaviour depends on the choices made of $\lambda,$ $\mu$ and $\rho.$ In fact, as explained in Section 3, the choice $\lambda=1000\kappa$ implies that $\alpha\simeq 0$ and the optimal square root uncertainty strategy approaches to the strategy solution of (\ref{3.17}), (\ref{3.18}) whose expected value for $\mu=0$ and $B>0$ is a convex function of time. As $p_0$ increases the parameter $B$ increases and the convexity of the optimal square root strategy increases. Otherwise, when $\lambda=1000\kappa$ the constant uncertainty strategy of Cheng et al. (2017)  approaches to the adaptive VWAP strategy whose expected value for $\mu=0$ is a linear function of time. \\In  Figure 2 we plot the sample trajectories of the optimal  square root uncertainty trading rate  (solid line) and of the optimal constant uncertainty trading rate (dotted line)  obtained with $p_0=10\%$ (left panel) and $p_0=30\%$ (right panel).  Looking at the sample trajectories of the optimal trading rates we can observe that  the optimal constant uncertainty rate is larger and more unstable than the optimal square root uncertainty rate and this effect is more evident towards the end of the liquidation interval where the optimal constant uncertainty rate spikes up significantly to achieve the full liquidation. This fact is expected because in the square root uncertainty it is possible to avoid the uncertainty choosing the trading rate equal to zero, otherwise in the constant uncertainty case this is not possible (see Bulthuis et al., 2017).

\section{Conclusions}
We have presented a new model of liquidation problem that takes into account execution risk. Under the assumption that execution risk affects both trading strategy and asset share price dynamics and that the magnitude of execution risk is proportional to the residual asset share position, we have modeled the liquidation problem as a linear quadratic stochastic optimal control problem and we have solved it. When the liquidation condition is enforced, i.e. the liquidation is completed at the final time of the liquidation interval, the optimal trading strategy is an ECIR square root process and belongs to the class of processes proposed by Delyon and Hu (2006); moreover the optimal trading rate found under execution risk is the optimal trading rate (without execution risk) of Almgren (2003) for a modified price.The model has the advantage of having explicit solution expressed by elementary functions obtained, differently from Cheng et al. (2017), without imposing any constraints on the parameters of the model.


\begin{thebibliography}{9}
	
	\bibitem{AC} Almgren, R., Chriss, N., Optimal execution of portfolio transactions, 2000, \textit{J. Risk}, 3(5), 5--39.
	
	\bibitem{A2003} Almgren, R., Optimal execution with nonlinear impact functions and trading enhanced risk, 2003,
	\textit{Appl. Math. Finance}, 10(1), 1--18.
	
	\bibitem{A2012} Almgren, R., Optimal trading with stochastic liquidity and volatility, \textit{SIAM J.  Financial Math.}, 2012, 3, 163--181.
	
	\bibitem{ABE} Ankirchner, S., Blanchet-Scalliet, C., Eyraud-Loisel, A., Optimal liquidation
	with additional information, \textit{Math. Financial Econ.}, Springer Verlag, 2016, 10(1), 1--14.
	
	
	\bibitem{BC} Bather, J.A.,  Chernoff, H., Sequential decisions in the control of a space ship (finite fuel),  \textit{J. Appl. Probab.}, 1967, 4(3), 584–-604.
	
	\bibitem{BSW} Bene\v{s}, V.E., Shepp, L.A., Witsenahusen, H.S., {Some solvable stochastic control problems}, \textit{Stoch.}, 1980, 4, 39--83.
	
	\bibitem{BCLW} Bulthuis, B., Concha, J., Leung, T., Ward, B., Optimal execution of limit and market orders with trade director, speed limiter, and fill uncertainty, \textit{Int. J. Financ. Eng.}, 2017, 4, 1--29.

	\bibitem{CDW2017} Cheng, X., Di Giacinto, M., Wang, T.-H.,  Optimal execution with uncertain order fills in Almgren–Chriss framework, \textit{Quant. Finance}, 2017, 17(1), 55--69.
	
	\bibitem{CDW2019}  Cheng, X., Di Giacinto, M., Wang, T.-H.,  Optimal execution with dynamic risk adjustment, \textit{J. Oper. Res. Soc.}, 2019, 70(10), 1662--1677.
	
	\bibitem{DH} Delyon, B., Hu, Y., Simulation of conditioned diffusion and application to parameter estimation, \textit{Stoch. Proc. Appl.}, 2006, 116, 1660--1675.
	
	\bibitem{DG} Durham, G.B., Gallant, A.R., Numerical techniques for maximum likelihood estimation of continuous-time diffusion processes, \textit{J. Bus. Econom. Statist.}, 2002, 20(3), 297--338.

       \bibitem{DZ} Du, S., Zhu, H., What is the optimal trading frequency in financial
	markets?, \textit{Rev. Econ. Stud.}, 2017, 84, 1606–-1651.
	
	\bibitem{FMRZ} Fatone, L., Mariani, F.,  Recchioni, M.C.,  Zirilli, F., A Trading execution model based on mean field games and optimal control, \textit{Appl. Math.}, 2014, 5(19), 3091--3116.
	
	\bibitem{FKTW} Forsyth, P.A., Kennedy, J.S., Tse, S.T., Windcliff, H., Optimal trade execution: A mean quadratic variation approach, \textit{J. Econ. Dyn. Control.}, 2012, 36(12), 1971--1991. 
	
	\bibitem{FW} Frei, C., Westray, N., Optimal Execution of a VWAP Order: a Stochastic Control Approach, \textit{Math. Finance}, 2015, 25(3), 612--639.
	
	\bibitem{GS}  Gatheral, J., Schied, A., Optimal trade execution under geometric brownian motion in the Almgren and Chriss    framework, \textit{Int. J. Theo. Appl. Finance}, 2011, 14(03), 53–368.
	
	\bibitem{GL} Gu\'{e}ant, O., Lehalle, C.A., General intensity shapes in optimal liquidation,  \textit{Math. Finance}, 2015, 25(3), 457--495.
	
	\bibitem{HW} Hull, J., White, A., Pricing Interest-Rate-Derivative Securities, \textit{Rev. Financial Studies}, 1990, 3(4), 573--392.
	
	\bibitem{Ka} Kalman, R.E., A new approach to linear filtering and prediction problems, \textit{J. Basic Eng.}, 1960, 82(1), 35--45.
	
	\bibitem{Kar} Karatzas, I., Probabilistic aspects of finite-fuel stochastic control \textit{Proc. Natl. Acad. Sci.}, USA, 1985, 82, 5579--5581.
	
	\bibitem{Kl} Klebaner, F.C., \textit{Introduction to Stochastic Calculus with Applications}, Imperial College Press, 2012. 
	
	\bibitem{KOW} Kyle, A.S., Obizhaeva,  A.A.,  Wang, Y., Smooth trading with overconfidence and market power, \textit{Rev.  Econ. Stud.}, 2017, 1, 1–-56.
	
	\bibitem{LS2013} Lorenz, C., Schied, A., Drift dependence of optimal trade execution strategies under transient price impact, \textit{Finance Stoch.}, 2012, 17(4), 743--777.
	

	\bibitem{SS} Sannikov, Y., Skrzypacz, A., Dynamic trading: price inertia and front-running,
	\textit{Stanford Univ. Grad. School of Business Research Paper}, No. 3487, 2016, 1--59.
	
	
	\bibitem{SST} Schied, A., Sch\"{o}neborn, T., Tehranchi, M., Optimal basket liquidation for CARA investors is deterministic, \textit{Appl. Math. Finance}, 2010, 17(6), 471--489.
	
	\bibitem{S2013} Schied, A., A control problem with fuel constraint and Dawson-Watanabe superprocesses, \textit{Ann. Appl. Probab.} , 2013, 23(6), 2472--2499.
	
	\bibitem{TFKW} Tse, S.T.,  Forsyth, P.A., Kennedy, J.S.,  Windcliff, H. , Comparison between the mean-variance optimal and the mean-quadratic-variation optimal trading strategies, \textit{Appl. Math. Finance}, 2013, 20(5), 415--449.
	
	\bibitem{Va} Vayanos, D., Strategic trading and welfare in a dynamic market, \textit{Rev. 
	Econ. Stud.}, 1999, 66 (2), 219--254.
	
	\bibitem{WGBS} Whitaker, G.A., Golightly, A., Boys, R.J., Sherlock, C., Improved bridge constructs for stochastic differential equations, \textit{Stat. Comput.}, 2017, 27, 885--900.
	
\end{thebibliography}
\end{document}